\documentclass{jams-l}
\usepackage{amssymb}
\usepackage{graphicx}

\newtheorem{theorem}{Theorem}[section]

\newtheorem{corollary}[theorem]{Corollary}
\newtheorem{proposition}[theorem]{Proposition}
\theoremstyle{definition}
\newtheorem{definition}[theorem]{Definition}
\newtheorem{example}[theorem]{Example}

\theoremstyle{remark}
\newtheorem{remark}[theorem]{Remark}

\numberwithin{equation}{section}

\usepackage{color}

\DeclareMathOperator{\Po}{Po}
\DeclareMathOperator{\Cl}{Cl}

\DeclareMathOperator{\Gal}{Gal}
\DeclareMathOperator{\Int}{Int}

\begin{document}

\title[P\'olya-Ostrowski Groups and Unit Indices in Real Biquadratic Fields]{
P\'olya-Ostrowski Groups and Unit Indices in Real Biquadratic Fields}

\author{Huda Naeem Hleeb Al-Jabbari}
\address{Department of Mathematics, Tarbiat Modares University, 14115-134, Tehran, Iran}
\curraddr{} \email{h-aljbbari@modares.ac.ir}
\thanks{}

\author{Abbas Maarefparvar$^{*}$}
\address{Department of Mathematics and Computer Science, University of Lethbridge, Lethbridge, Canada}
\curraddr{}
\email{abbas.maarefparvar@uleth.ca}
\thanks{$^{*}$Corresponding author}
\date{}
\dedicatory{}
\commby{}

%\author{Abbas Maarefparvar}
%\address{School of Mathematics, Institute for Research in Fundamental Sciences (IPM), P.O. Box: 19395-5746, Tehran, Iran.}
%\curraddr{}
 %\email{a.marefparvar@ipm.ir}
%\thanks{This research was supported by a grant from IPM}
%\date{}
%\dedicatory{}
\commby{}
%\author{Ali Rajaei$^{**}$}
%\address{Department of Mathematics, Tarbiat Modares University, 14115-134, Tehran, Iran}
%\curraddr{}
%\email{alirajaei@modares.ac.ir}
%\thanks{$^{**}$Corresponding author}
%\author{Ehsan Shahoseini}
%\address{Department of Mathematics, Tarbiat Modares University, 14115-134, Tehran, Iran}
%\curraddr{} \email{ehsan\_shahoseini@modares.ac.ir}
%\thanks{}
\maketitle

%    \subjclass is required.
%\subjclass[2010]{Primary 11R04, 11R16, 11R29, 11R34, 11R37}

%\date{}

%\dedicatory{}

%    "Communicated by" -- provide editor's name; required.
%%\commby{}

%    Abstract is required.
%%%\maketitle
%\begin{abstract}
{\noindent \bf{Abstract:}}
The P\'olya-Ostrowski group of a Galois number field $K$, is the subgroup $\Po(K)$ of the ideal class group $\Cl(K)$ of $K$ generated by the classes of  all the strongly ambiguous ideals of $K$. The number field $K$ is called a P\'olya field, whenever $\Po(K)$ is trivial. In this paper, using  some results of Bennett Setzer \cite{Bennett} and Zantema \cite{Zantema}, we  give an explicit relation between the order of P\'olya groups and the Hasse unit indices in real biquadratic fields. 
As an application, we refine Zantema's upper bound on the number of ramified primes in  P\'olya real biquadratic fields.

\vspace{.2cm} {\noindent \bf{Keywords:}}~ P\'olya fields, P\'olya groups, Biquadratic fields, Hasse unit index, Galois cohomology.

\vspace{.2cm} {\noindent \bf{Notation.}}~
For a number field $K$, we denote by $d_K$, $\Cl(K)$, $h_K$, $\mathcal{O}_K$ and $U_K$ the discriminant of $K$, the ideal class group, the class number, the ring of integers and  the group of units of $K$, respectively. For a prime number $p$, we denote by $e_p(K/\mathbb{Q})$ the ramification index of $p$ in $K/\mathbb{Q}$. The notation $N_{K/\mathbb{Q}}$ denotes both  the element and ideal norm morphisms from $K$ to $\mathbb{Q}$.
%For a finite extension $L/K$ of number fields, we denote by $\mathcal{N}_{L/K}$  both the element and ideal norm morphisms from $L$ to $K$.

\section{Introduction}\label{intro}

The study of integer-valued polynomials on rings of integers of number fields
originates in the seminal works of P\'olya \cite{Polya} and Ostrowski \cite{Ostro}.
For a number field $K$ with ring of integers $\mathcal{O}_K$, let
\[
\mathrm{Int}(\mathcal{O}_K)
=
\{ f\in K[X] : f(\mathcal{O}_K)\subseteq \mathcal{O}_K \}
\]
denote the ring of integer-valued polynomials on $\mathcal{O}_K$.
It is known that $\mathrm{Int}(\mathcal{O}_K)$ is a free $\mathcal{O}_K$-module
\cite[Section~2]{Zantema}.
Following Zantema \cite{Zantema}, the field $K$ is called a \emph{P\'olya field} if
$\mathrm{Int}(\mathcal{O}_K)$ admits a regular basis, that is, a basis containing
exactly one polynomial of each degree. An equivalent arithmetic formulation is given in terms of the \emph{Ostrowski ideals}.
For each positive integer $q$, define
\[
\Pi_q(K)
=
\prod_{\substack{\mathfrak{m}\in\mathrm{Max}(\mathcal{O}_K)\\
		N_{K/\mathbb{Q}}(\mathfrak{m})=q}}
\mathfrak{m},
\]
with the convention that $\Pi_q(K)=\mathcal{O}_K$ if no such ideals exist. The field $K$ is P\'olya if and only if all ideals $\Pi_q(K)$ are principal.

Cahen and Chabert \cite{Cahen-Chabert's book} introduced the \emph{P\'olya--Ostrowski group}
(or simply the \emph{P\'olya group}) $\mathrm{Po}(K)$ as the subgroup of the ideal
class group $\mathrm{Cl}(K)$ generated by the classes of all Ostrowski ideals.
Thus $K$ is P\'olya if and only if $\mathrm{Po}(K)$ is trivial.

When $K/\mathbb{Q}$ is Galois, the P\'olya group admits a cohomological description.
In this case, the Ostrowski ideals generate precisely the subgroup of strongly
ambiguous ideals, and $\mathrm{Po}(K)$ coincides with the subgroup of
$\mathrm{Cl}(K)$ consisting of strongly ambiguous ideal classes. The reader is referred to \cite{ChabertI,  Leriche 2013, MaarefJNT2021} for some results on P\'olya fields and P\'olya groups.

Zantema \cite{Zantema} established an exact sequence relating $\mathrm{Po}(K)$
to the cohomology group $H^{1}(\mathrm{Gal}(K/\mathbb{Q}),U_K)$ and to the
ramification indices of primes in $K$.

\begin{proposition}\cite[Proposition 3.1]{Zantema} \label{proposition, Zantema's exact sequence}
	For a Galois extension $K/\mathbb{Q}$, the following sequence is exact:
	\begin{equation} \label{equation, Zantema exact sequence}
	0 \longrightarrow H^1(\Gal(K/\mathbb{Q}),U_K) \longrightarrow \bigoplus _{\text{p prime}} \mathbb{Z} / e_p(K/\mathbb{Q}) \mathbb{Z} \longrightarrow \Po(K) \longrightarrow 0,
	\end{equation}
	where $e_p(K/\mathbb{Q})$ denotes the ramification index of a prime $p$ in $K$.
\end{proposition}

The above exact sequence allows one to derive restrictions on ramification in P\'olya fields. In particular, Zantema \cite[Section 5]{Zantema} showed that in a P\'olya real biquadratic field
at most five primes ramify, and that this bound is sharp.
For instance $K=\mathbb{Q}(\sqrt{35}, \sqrt{381})$ is P\'olya with five ramifications \cite[page 20]{Zantema} (Some infinite families of P\'olya real biquadratic fields with maximum ramification have been found in \cite{MaarefJNT2021}). For a real biquadratic field $K$, Bennett Setzer \cite{Bennett} described the structure of $H^1(\Gal(K/\mathbb{Q}),U_K)$ in terms of unit-theoretic data from the quadratic subfields $k_1,k_2,k_3$, a description that has been used in several subsequent works, including earlier paper of the second author \cite{MaarefJNT2021} The present paper builds on these foundational results but differs in scope and outcome: rather than focusing on particular families of fields, we obtain an explicit and uniform formula for the order of the P\'olya group of an arbitrary real biquadratic fields $K$ in terms of the
\emph{Hasse unit index} $	\left[U_K : U_{k_1} U_{k_2} U_{k_3}\right]$, ramification indices, and the number of quadratic subfields whose fundamental units have negative norm, see Theorem \ref{theorem, relatione between order of Polya group and unit index}. This leads, as an application, to refined upper bounds on the number of ramified primes in Pólya real biquadratic fields, improving Zantema’s classical bound in several cases, see Corollary \ref{corollary, improving Zantema upper bound for Polya biquadratic fields}.

\section{An Explicit Relation Between the Order of the P\'olya Groups and the Hasse Unit Indices of Real Biqauadratic Fields} \label{Setzer's results}
In this section, summarizing some results of Bennett Setzer \cite{Bennett} for real biqudratic fields $K$, we find an explicit relation between $\# \Po(K)$ and the Hasse unit index of $K$.

\begin{proposition}\label{proposition, Setzer results}\cite{Bennett}
%Let $d_1,d_2 (\neq 1)$ be two different squarefree positive integers and $K=\mathbb{Q}(\sqrt{d_1},\sqrt{d_2})$
Let $K=\mathbb{Q}(\sqrt{d_1},\sqrt{d_2})$ be a real biquadratic field with Galois group $G$, where $d_1,d_2 \neq 1$ are two different square-free positive integers.
Denote by $k_1,k_2,k_3$ the three real quadratic subfields of $K$. Let  $H=H^1(G,U_K)$, and $H[2]$ be the $2$-torsion subgroup of $H$. Then,
%Let $d_1,d_2 (\neq 1)$ be two different squarefree positive integers and
%$K=\mathbb{Q}(\sqrt{d_1},\sqrt{d_2})$. Identify $\Gal(K/\mathbb{Q})$
% as $G=\{1,S_1,S_2,S_3\}$ where $k_i:=\mathbb{Q}(\sqrt{d_i})$ is the fix field of $S_i$ ($d_3$ is the square-free part of $d_1 d_2$).
% Let  $H=H^1(G,U_K)$, $H[2]$ be the $2$-torsion subgroup of $H$, and
%{\small
%\begin{equation*}
%\mathcal{G}=\{ (u_1,u_2,u_3) \in U_{k_1} \times U_{k_2} \times U_{k_3} \mid N_{k_1/\mathbb{Q}}(u_1)=N_{k_2/\mathbb{Q}}(u_2)=N_{k_3/\mathbb{Q}}(u_3)=\text{sign}(u_1u_2u_3)\}.
%\end{equation*}}
the following assertions hold.
\begin{itemize}
\item[(i)] \cite[Theorem 4]{Bennett} $\left[H:H[2] \right] \leq 2$. Moreover $\left[H:H[2] \right]=2$ if and only if $2$ ramifies totally in $K/\mathbb{Q}$ and all $k_i$'s contain elements with the same norm either $2$ or $-2$.

\item[(ii)]\cite[$\S$ 4,5]{Bennett} Denote by $t$ the number of quadratic subfields of $K$ whose fundamental units have negative norm. Then
\begin{equation*}
  	\left[U_K : U_{k_1} U_{k_2} U_{k_3}\right]. \#H[2]  = \begin{cases}
               2^{5-t}              & : \,  t =0,1, \\ 
              2^{3}            & : \,  t=2,3.
           \end{cases}
\end{equation*}
Moreover for $t=3$, $H[2]$ contains a subgroup of order $4$.

\item[(iii)] For $i=1,2,3$, let $\alpha_i$ be the generator of units in $k_i$ with positive norms. If at least one $k_i$ has no units of negative norm, then $H[2]$ is isomorphic to the subgroup
\begin{equation} \label{equation, subgroup Q*}
\mathcal{A}	=\langle [d_1],[d_2],[N_{k_1/\mathbb{Q}}(\alpha_1+1)],[N_{k_2/\mathbb{Q}}(\alpha_2+1)], [N_{k_1/\mathbb{Q}}(\alpha_3+1)] \rangle
\end{equation}
of the quotient group $\mathbb{Q}^{\times}/\mathbb{Q}^{\times 2}$. Otherwise, $H[2]$ would be isomorphic to the subgroup $\mathcal{A}$ with possibly one more generator.
\end{itemize}
\end{proposition}

For a real biquadratic field $K$,  the first cohomology group $H^1(\Gal(K/\mathbb{Q}),U_K)$ can be seen as a ``common term'' in the above results of Bennet Setzer and the exact sequence \eqref{equation, Zantema exact sequence}. This observation leads us to relate the order of P\'olya groups of real biquadratic fields to their Hasse unit indices.

\begin{theorem} \label{theorem, relatione between order of Polya group and unit index}
Let $K=\mathbb{Q}(\sqrt{d_1},\sqrt{d_2})$ be a real biquadratic field with quadratic subfields $k_i=\mathbb{Q}(\sqrt{d_i})$ ($d_3$ is the squarefree part of $d_1 d_2$). Denote by $t$ the number of quadratic subfields of $K$ whose fundamental unit have negative norm. 
\begin{itemize}
\item[(i)] If $2$ ramifies totally in $K/\mathbb{Q}$ and all $k_i$'s contain elements with the same norm either $2$ or $-2$, then:
\begin{equation*}
    \# \Po(K) = \begin{cases}
               \frac{	\left[U_K : U_{k_1} U_{k_2} U_{k_3}\right]. \prod_{p \mid d_K}e_p(K/\mathbb{Q})}{2^{6-t}}              & : \,  t = 0,1, \\
                & \\ 
               \frac{	\left[U_K : U_{k_1} U_{k_2} U_{k_3}\right]. \prod_{p \mid d_K}e_p(K/\mathbb{Q})}{2^4}              & : \, t =2,3.\\
           \end{cases}
\end{equation*}

\item[(ii)] Otherwise:
\begin{equation*}
    \# \Po(K) = \begin{cases}
               \frac{	\left[U_K : U_{k_1} U_{k_2} U_{k_3}\right]. \prod_{p \mid d_K}e_p(K/\mathbb{Q})}{2^{5-t}}              & : \,  t = 0,1, \\
               & \\
               \frac{	\left[U_K : U_{k_1} U_{k_2} U_{k_3}\right]. \prod_{p \mid d_K}e_p(K/\mathbb{Q})}{2^3}              & : \, t=2,3.\\
           \end{cases}
\end{equation*}
\end{itemize}
In particular, for $t\in \{2,3\}$  P\'olya groups have the same order.
\end{theorem}

\begin{proof}
We only prove part (i), and the second part is proved similarly. Suppose that $2$ ramifies totally in $K/\mathbb{Q}$, and all $k_i$'s contain elements with same norm either $2$ or $-2$. Then, by the first two parts of Proposition \ref{proposition, Setzer results}, we have
\begin{equation} \label{eq1}
\# H^1(\Gal(K/\mathbb{Q}),U_K) \cdot \left[ U_K : U_{k_1} U_{k_2} U_{k_3}\right]=
 \begin{cases}
		2^{6-t}              & : \,  t =0,1, \\ 
		2^{4}            & : \,  t=2,3.
	\end{cases}
\end{equation}
On the other hand, by the exact sequence \eqref{equation, Zantema exact sequence}, 
\begin{equation} \label{eq2}
	\# \Po(K)=\frac{\prod_{p \mid d_K}e_p(K/\mathbb{Q})}{\# H^1(\Gal(K/\mathbb{Q}),U_K)}.
\end{equation}
Using \eqref{eq1} with \eqref{eq2}, we get the assertion.

%Combining parts $(i)-(ii)$ in Proposition \eqref{proposition, Setzer results} with
%Zantema's exact sequence in Proposition \eqref{proposition, Zantema's exact sequence}, we get the assertions.
\end{proof}

For a real biquadratic field $K$, we can use Proposition \eqref{proposition, Setzer results} to find a lower bound for the order of the first cohomology group $H^1(\Gal(K/\mathbb{Q}),U_K)$. Then using the exact sequence \eqref{equation, Zantema exact sequence} we can compute (or at least approximate) the order of $\Po(K)$, and then Theorem \eqref{theorem, relatione between order of Polya group and unit index} gives us  the corresponding Hasse unit index of $K$. To do so, we also need the following property of a generator of positive norm units in real quadratic fields.

\begin{proposition} \cite{MaarefJNT2021} \label{proposition, properties of m_k}
	Let $k=\mathbb{Q}(\sqrt{d})$ be a real quadratic field and $\alpha_k >1$ be the generator of positive norm units in $k$. Assume that $n_k$ is the squarefree part of $N_{k/\mathbb{Q}}(\alpha_{k}+1)$. Then $n_k \not \in \{ 1,d\}$, and $n_k$ divides the discriminant of $k$.
	%, $n_k$ is the norm of an integer in $k$, and $n_k .\epsilon$ is a square in $k$.
\end{proposition}

%For a real biquadratic field $K$ with Galois group $G$,
%Propositions \eqref{proposition, properties of m_k} and \eqref{proposition, Setzer results} and Corollary \eqref{corollary, isomorphism between H[2] and rho(G)} introduce a practical method for estimation $\#H^1(G,U_K)$.  Then using the exact sequence \eqref{equation, Zantema exact sequence} we can compute (or at least approximate) the order of $\Po(K)$, and then Theorem \eqref{theorem, relatione between order of Polya group and unit index} gives us  the corresponding unit index of $K$.

\begin{example}
	Let  $p \equiv q \equiv 1 (\mathrm{mod}\, 4)$ be distinct  prime numbers, $K=\mathbb{Q}(\sqrt{p},\sqrt{q})$,  $k_1=\mathbb{Q}(\sqrt{p})$, $k_2=\mathbb{Q}(\sqrt{q})$, $k_3=\mathbb{Q}(\sqrt{pq})$ and denote the fundamental unit of $k_i$ by $\epsilon_i$, respectively. Then $N_{k_1/\mathbb{Q}}(\epsilon_1)=N_{k_2/\mathbb{Q}}(\epsilon_2)=-1$, see e.g. \cite[Proposition 4]{Kis}. By Theorem \eqref{theorem, relatione between order of Polya group and unit index}, we get
	\begin{equation} \label{equation, Po(K) is 2 times of Unit index}
	2. \# \Po(K)= \left[U_K : U_{k_1} U_{k_2} U_{k_3}\right].
	\end{equation}

	Let $H=H^1(\Gal(K/\mathbb{Q}),U_K)$. Since $2$ doesn't ramify totally in $K/\mathbb{Q}$, part $(i)$ in Proposition \eqref{proposition, Setzer results} implies that $H=H[2]$. Using Propositions \ref{proposition, Zantema's exact sequence}, \ref{proposition, Setzer results} and \ref{proposition, properties of m_k}, we have
	\begin{itemize}
		\item[•] if $N_{k_3/\mathbb{Q}}(\epsilon_3)=+1$, then
		% and Corollary \eqref{corollary, isomorphism between H[2] and rho(G)}:
		\begin{equation*}
			\left< p,q  \right> \, (\mathrm{mod}\, \mathbb{Q}^{\times 2})  \subseteq \mathcal{A} \simeq H[2] = H \subseteq \frac{\mathbb{Z}}{2\mathbb{Z}}\oplus \frac{\mathbb{Z}}{2\mathbb{Z}};
		\end{equation*}
		\item[•] if $N_{k_3/\mathbb{Q}}(\epsilon_3)=-1$, then
		\begin{equation*}
			\frac{\mathbb{Z}}{2\mathbb{Z}}\oplus \frac{\mathbb{Z}}{2\mathbb{Z}} \subseteq H[2] =H  \subseteq \frac{\mathbb{Z}}{2\mathbb{Z}}\oplus \frac{\mathbb{Z}}{2\mathbb{Z}}.
		\end{equation*}
	\end{itemize}
	
	Hence, in both cases, $ H^1(\Gal(K/\mathbb{Q}),U_K) \simeq \frac{\mathbb{Z}}{2\mathbb{Z}}\oplus \frac{\mathbb{Z}}{2\mathbb{Z}}$ and by  the exact sequence \eqref{equation, Zantema exact sequence}, $K$ is a P\'olya field. By the relation  \eqref{equation, Po(K) is 2 times of Unit index}, we find $\left[U_K : U_{k_1} U_{k_2} U_{k_3}\right]=2$.
\end{example}

\begin{remark}
 Note that in the above example $k_1$ and $k_2$ are P\'olya. Indeed, in \textit{almost} cases, compositum of two  quadratic P\'olya fields is a P\'olya field  \cite[$\S$ 5]{Leriche 2013} for which unit index is obtained immediately from
  Theorem \eqref{theorem, relatione between order of Polya group and unit index}.
\end{remark}

As mentioned before, Zantema \cite{Zantema} gave a sharp upper bound on the number of ramified primes in P\'olya real biquadratic fields.

\begin{proposition} \cite[Section 5]{Zantema}
Let $K$ be a real biquadratic field, and denote by $s_K$ the number of primes that are ramified in $K/\mathbb{Q}$. If $K$ is P\'olya, then $s_K \leq 5$. Moreover, this upper bound is sharp. 
\end{proposition}

As a direct consequence of Theorem \ref{theorem, relatione between order of Polya group and unit index}, we obtain a \emph{refined} version of Zantem's upper bound in the above proposition.

%As an application of Theorem \ref{theorem, relatione between order of Polya group and unit index}, we can improve Zantema's upper bound on the number of ramified primes in P\'olya real biquadrtaic fields.
%A direct consequence of Theore

%Theorem \eqref{theorem, relatione between order of Polya group and unit index} also gives an upper bound on the number of ramified primes in P\'olya real biquadratic fields that can be seen as an improvement on Zantema's upper bound \cite[$\S$5]{Zantema}, see Remark \eqref{remark, Zantema upper bound for Polya biquadratic fields}:

\begin{corollary} \label{corollary, improving Zantema upper bound for Polya biquadratic fields}
	Let
$K$ be a real biquadratic field. Denot by $t$ the number of quadratic subfields of $K$ whose fundamental units have negative norm, and by $s_K$ the number of primes that are ramified in $K/\mathbb{Q}$. If $K$ is P\'olya, then
\begin{equation*}
s_K \leq \left\{
\begin{array}{rl}
5 \, :& \, t=0, \\
 4 \, :& \,   t=1, \\
 3 \, :& \,  t=2,3.
\end{array} \right.
\end{equation*}
\end{corollary}

\begin{proof}
 Suppose that $\Po(K)=0$. We prove the cases $t=2,3$. The other two cases, say $t=0,1$, can be proved similarly. Let $k_1,k_2,k_3$ be the three quadratic subfield of $K$, and suppose that $k_1$ and $k_2$ have some units of negative norm. If $2$ is totally ramified in $K/\mathbb{Q}$, and all $k_i$'s contain elements with the same norm either $2$ or $-2$, then by part (i) of Theorem \ref{theorem, relatione between order of Polya group and unit index}, we get
\begin{equation*}
\prod_{p | d_K} e_p(K/\mathbb{Q})=e_2(K/\mathbb{Q}) \cdot \prod_{\substack{p | d_K\\ p \neq 2}} e_p(K/\mathbb{Q}) = \frac{2^4}{\left[U_K : U_{k_1} U_{k_2} U_{k_3}\right]},
\end{equation*}
which implies that 
\begin{equation*}
 \prod_{\substack{p | d_K\\ p \neq 2}} e_p(K/\mathbb{Q}) \leq 2^2,
\end{equation*}
or equivalently, $s_K \leq 3$. On the other hand, if $2$ doesn't ramify totally in $K/\mathbb{Q}$ or one of the equations
\begin{equation*}
N_{k_i/\mathbb{Q}}(x_i)=+2 \quad \text{or} \quad -2, \quad i=1,2,3,
\end{equation*}
doesn't have any solution $x_i \in \mathcal{O}_{k_i}$, then by part (ii) of Theorem \ref{theorem, relatione between order of Polya group and unit index}, we obtain
\begin{equation*}
	\prod_{p | d_K} e_p(K/\mathbb{Q}) = \frac{2^3}{\left[U_K : U_{k_1} U_{k_2} U_{k_3}\right]} \leq 2^3,
\end{equation*}
and again we get $s_K \leq 3$.
\end{proof}

%\section*{Acknowledgment} 
%The first author would like to thank Jonas Morrissey for a translation of Hayashi's paper and Abbas Maarefparvar for a careful reading of a first draft of this paper. The second author's research
%was in part supported by a grant from IPM (No. 96110120).


\bibliographystyle{amsplain}

\end{document}